\documentclass[a4paper]{elsarticle}
\usepackage{multirow}
\usepackage[fleqn]{amsmath}
\usepackage{graphicx}
\usepackage{caption}
\usepackage{dsfont}
\usepackage{color}
\usepackage{amsfonts}
\usepackage{subfigure}
\usepackage{lineno,hyperref}
\modulolinenumbers[5]
\usepackage{CJK}
\usepackage{indentfirst}
\usepackage{amsmath}
\usepackage{geometry}

\journal{\textbf{Computers and Mathematics with Applications}}

\begin{document}

\begin{frontmatter}


\title{Stochastic Burgers equation with fractional derivative driven by multiplicative noise}

\author[mymainaddress]{Guang-an Zou\corref{mycorrespondingauthor}}
\cortext[mycorrespondingauthor]{Corresponding author}
\ead{zouguangan00@163.com}

\author[mymainaddress]{Bo Wang}

%\author{Guang-an Zou\fnref{myfootnote}}
%\address{School of Mathematics and Statistics, Henan University, Kaifeng 475004, P. R. China}

\address[mymainaddress]{School of Mathematics and Statistics, Henan University, Kaifeng 475004, China}

\begin{abstract}

This article is devoted to the study of the existence and uniqueness of mild solution to time- and space-fractional stochastic Burgers equation perturbed by multiplicative white noise. The required results are obtained by stochastic analysis techniques, fractional calculus and semigroup theory. We also proved the regularity properties of mild solution for this generalized Burgers equation.

\end{abstract}

\begin{keyword}
Stochastic Burgers equation, fractional derivative, mild solution, regularity properties.
%\texttt{elsarticle.cls}\sep \LaTeX\sep Elsevier \sep template
%\MSC[2010] 00-01\sep  99-00
\end{keyword}

\end{frontmatter}

%\linenumbers

\section{Introduction}

Stochastic Burgers equation (SBE) plays an important role in the modeling of many phenomena in different fields, such as fluid dynamics, nonlinear acoustics, hydrodynamics, cosmology, astrophysics and statistical physics, and so on. In the last decade, SBE has gained a great development in both theory and application and a large volume of literature is available on this subject (see e.g.[1-4] and references therein). It is particularly mentioned that when the Laplacian operator $\Delta$ in SBE is replaced by fractional derivative, which can be used to describe anomalous diffusion processes in fractal flow and acoustic waves propagation in porous media [6,22,23]. Sugimoto [5] have studied the generalized Burgers-type equation with a fractal power of Laplacian in the principal part, which described the unidirectional propagation of acoustic waves through a gas-filled tube with a boundary layer. Besides, the space-fractional SBE also can be used to study the acoustic waves propagation in tunnels during the passage of the trains, which may yield a memory effect and other types of resonance phenomena [6]. On the other hand, time-fractional differential equations are found to be quite effective in modelling anomalous diffusion processes as its can characterize the long memory processes [14-16,21,24]. Hence, Burgers equation with time-fractional can be adapted to describe the memory effect of the wall friction through the boundary layer [7]. Furthermore, the analytical solutions of the time- and space-fractional Burgers equations have been investigated by variational iteration method [7] and Adomian decomposition method [8].

In this study, we focus on the following generalized SBE with time-space fractional derivative on a bounded domain $D\subset \mathds{R}^{d}(1\leq d\leq3)$:

\begin{align*}
 ^{C}D_{t}^{\beta}u-u\cdot\nabla u+(-\triangle)^{\frac{\alpha}{2}}u=g(u)\dot{W}(t),(t,x)\in(0,T]\times D, \tag{1.1}
\end{align*}
subject to the initial condition:
\begin{align*}
 u(0,x)=u_{0}(x),x\in D, \tag{1.2}
\end{align*}
and the Dirichlet boundary conditions:
\begin{align*}
 u(t,x)|_{\partial D}=0,t\in[0,T], \tag{1.3}
\end{align*}   
in which the term $g(u)\dot{W}(t)=g(u)\frac{d W(t)}{d t}$ describes a state dependent random noise, where $W(t)_{t\in[0,T]}$ is a $\mathcal{F}_{t}$-adapted Wiener process defined on a completed probability space $(\Omega,\mathcal{F},\mathds{P})$ with the expectation $\mathbb{E}$, and associate with the normal filtration $\mathcal{F}_{t}=\sigma\{W(s):0\leq s\leq t\}$; The operator $(-\triangle)^{\frac{\alpha}{2}},\alpha\in(1,2)$ stands for the fractional power of the Laplacian (see [25]); We denote by $^{C}D_{t}^{\beta}$ the Caputo derivative of order $\beta$, which is defined by (see [9])
\begin{align*}
^{C}D_{t}^{\beta}u(t,x)=\begin{cases}
\frac{1}{\Gamma(1-\beta)}\int_{0}^{t}\frac{\partial u(s,x)}{\partial s}\frac{ds}{(t-s)^{\beta}},~0<\beta<1,\\
\frac{\partial u(t,x)}{\partial t}, \qquad \qquad \qquad \qquad \beta =1,\\
\end{cases} \tag{1.4}
\end{align*}
where $\Gamma(\cdot)$ stands for the gamma function $\Gamma(\beta)=\int_{0}^{\infty}t^{\beta-1}e^{-t}dt$.

Eq.(1.1) might be used to model anomalous diffusion processes in disordered media, and describe the acoustic wave propagation in porous media with memory effect and with random effects. Notice that the study of space-fractional SBE can be found in some literatures. For details, Brze\'{z}niak and Debbi [10] proved existence and uniqueness of a mild global solution to the Cauchy problem for the stochastic fractional Burgers equation. Brze\'{z}niak et al.[6] studied the ergodic properties of the solution for space-fractional SBE. Yang [11] proposed some estimates on the solution of space-fractional SBE and given the invariant measure. Lv and Duan [12] investigated the existence of martingale solutions and weak solutions for  space-fractional SBE on a bounded domain. However, to the best of our knowledge, there are no existing works for the time- and space-fractional SBE, which is a fascinating and useful problem.

The main contribution of this paper is to establish the existence, uniqueness, and regularity properties of mild solution to time-space fractional SBE driven by multiplicative noise, which generalizes many previous works [6,11,12]. The rest of the paper is organized as follows. In Section 2, we will introduce some notations and preliminaries, which play a crucial role in our theorem analysis. In Section 3, the existence and uniqueness of mild solution to the problems of time-space fractional SBE are obtained by stochastic analysis techniques, fractional calculus and semigroup theory. Finally, the spatial and temporal regularity properties of mild solution to this time-space fractional SBE are proved.

\section{Notations and preliminaries}

Denote the basic functional space $L^{p}(D),1\leq p<\infty$ and $H^{s}(D)$ by the usual Lebesgue and Sobolev spaces, respectively. We assume that $A$ is the negative Laplacian $-\Delta$ in a bounded domain $D$ with zero Dirichlet boundary conditions in a Hilbert space $H=L^{2}(D)$, which are given by
\begin{align*}
A=-\Delta,~\mathcal{D}(A)=H_{0}^{1}(D)\cap H^{2}(D).
\end{align*}

Since the operator $A$ is self-adjoint on $H$ with discrete spectral, i.e., there exists the eigenvectors $e_{n}$ with corresponding eigenvalues $\lambda_{n}$ such that
\begin{align*}
Ae_{n}=\lambda_{n}e_{n},e_{n}=\sqrt{2}\sin(n\pi),\lambda_{n}=\pi^{2}n^{2},n\in N^{+}.
\end{align*}

For any $s>0$, let $\dot{H}^{s}$ be the domain of the fractional power $A^{\frac{s}{2}}=(-\triangle)^{\frac{s}{2}}$, which can be defined by
\begin{align*}
A^{\frac{s}{2}}e_{n}=\lambda_{n}^{\frac{s}{2}}e_{n}, n=1,2,\ldots,
\end{align*}
and
\begin{align*}
 \dot{H}^{s}=\mathcal{D}(A^{\frac{s}{2}})=\{v\in L^{2}(D),s.t.\|v\|_{\dot{H}^{s}}^{2}=\sum\limits_{n=1}^{\infty}\lambda_{n}^{\frac{s}{2}}v_{n}^{2}<\infty\},
\end{align*}
where $v_{n}:=\langle v,e_{n}\rangle$ with the inner product $\langle \cdot,\cdot\rangle$ in $L^{2}(D)$. We denote that $\|v\|_{\dot{H}^{s}}=\|A^{\frac{s}{2}}v\|$, and the corresponding dual space $\dot{H}^{-s}$ with the inverse operator $A^{-\frac{s}{2}}$. We also denote $A_{s}$ for $A^{\frac{s}{2}}$ and the bilinear operator $B(u,v)=u\cdot\nabla v$, and $\mathcal{D}(B)=H_{0}^{1}(D)$ with a slight abuse of notation $B(u):=B(u,u)$. Then the Eqs.(1.1)-(1.3) can be rewritten as the following abstract formulation:
\begin{align*}
\begin{cases}
^{C}D_{t}^{\beta}u(t)=-A_{\alpha}u(t)+B(u(t))+g(u(t))\frac{d W(t)}{d t},t>0,\\
u(0)=u_{0},
\end{cases} \tag{2.1}
\end{align*}
where $\{W(t)\}_{t\geq0}$ is a $Q$-Wiener process with linear bounded covariance operator $Q$ such that $\mathrm{Tr} (Q)<\infty$. Further, there exists the eigenvalues $\lambda_{n}$ and corresponding eigenfunctions $e_{n}$ satisfy $Q e_{n}=\lambda_{n}e_{n},n=1,2,\ldots$, then the Wiener process is given by
\begin{align*}
W(t)=\sum\limits_{n=1}^{\infty}\lambda^{1/2}_{n}\beta_{n}(t)e_{n},
\end{align*}
in which $\{\beta_{n}\}_{n\geq1}$ is a sequence of real-valued standard Brownian motions. 

Let $L_{0}^{2}=L^{2}(Q^{1/2}(H),H)$ denote the space of Hillbert-Schmidt operators from $Q^{1/2}(H)$ to $H$ with the norm $\|\phi\|_{L_{0}^{2}}:=\|\phi Q^{1/2}\|_{HS}=(\sum\limits_{n=1}^{\infty}\phi Q^{1/2}e_{n})^{1/2}$, i.e., $L_{0}^{2}=\{\phi\in L(H):\sum\limits_{n=1}^{\infty}\|\phi Q^{1/2}e_{n})\|^{2}<\infty\}$, where $L(H)$ is the space of bounded linear operators from $H$ to $H$.

For an arbitrary Banach space $B$, we denote $\|\cdot\|_{L^{p}(\Omega;B)}$ by the norm in $L^{p}(\Omega,\mathcal{F},\mathds{P};B)$, which defined as
\begin{align*}
\|v\|_{L^{p}(\Omega;B)}=(\mathbb{E}[\|v\|_{B}^{p}])^{\frac{1}{p}},\forall ~v\in L^{p}(\Omega,\mathcal{F},\mathds{P};B),
\end{align*}
for any $p\geq2$.

We shall also need the following result with respect to the fractional operator $A_{\alpha}$ (see Ref.[11]).

\textbf{Lemma 2.1.} For any $\alpha>0$, an analytic semigroup $S_{\alpha}(t)=e^{-tA_{\alpha}},t\geq0$ is generated by the operator $-A_{\alpha}$ on $L^{p}$, and for any $\nu\geq0$, there exists a constant $C_{\alpha,\nu}$ dependent on $\alpha$ and $\nu$ such that
\begin{align*}
\|A_{\nu}S_{\alpha}(t)\|_{\mathcal{L}(L^{p})}\leq C_{\alpha,\nu}t^{-\frac{\nu}{\alpha}}, t>0,\tag{2.2}
\end{align*}
in which $\mathcal{L}(B)$ denotes the Banach space of all linear bounded operators from $B$ to itself.

Next, we will introduce the following lemma to estimate the stochastic integrals, which contains the Burkholder-Davis-Gundy's inequality.

\textbf{Lemma 2.2.}([13]) For any $0\leq t_{1}<t_{2}\leq T$ and $p\geq2$, and for any predictable stochastic process $v:[0,T]\times\Omega\rightarrow L_{0}^{2}$, which satisfies
\begin{align*}
\mathbb{E}[(\int_{0}^{T}\|v(s)\|_{L_{0}^{2}}^{2}ds)^{\frac{p}{2}}]<\infty,
\end{align*}
then we have
\begin{align*}
\mathbb{E}[\|\int_{t_{1}}^{t_{2}}v(s)dW(s)\|^{p}]\leq C(p)\mathbb{E}[(\int_{t_{1}}^{t_{2}}\|v(s)\|_{L_{0}^{2}}^{2}ds)^{\frac{p}{2}}], \tag{2.3}
\end{align*}
where $C(p)=[\frac{p(p-1)}{2}]^{\frac{p}{2}}(\frac{p}{p-1})^{p(\frac{p}{2}-1)}$ is a constant.

Inspired by the definition of the mild solution to the time-fractional differential equations (see Refs.[14-18]), we give the following definition of mild solution for our time-space fractional stochastic Burgers equation.

\textbf{Definition 2.1.}  A $\mathcal{F}_{t}$-adapted process $(u(t))_{t\in[0,T]}$ is called a mild solution to (2.1), if $(u(t))_{t\in[0,T]}\in C([0,T];\dot{H}^{\nu})$ $\mathds{P}$-a.e., and
it holds
\begin{align*}
 u(t)&=E_{\beta}(t)u_{0}+\int_{0}^{t}(t-s)^{\beta-1}E_{\beta,\beta}(t-s)B(u(s))ds\\
 &\hspace{2mm}+\int_{0}^{t}(t-s)^{\beta-1}E_{\beta,\beta}(t-s)g(u(s))dW(s), \tag{2.4}
\end{align*}
for a.s. $\omega\in\Omega$, where the generalized Mittag-Leffler operators $E_{\beta}(t)$ and $E_{\beta,\beta}(t)$ are defined as
\begin{align*}
E_{\beta}(t)=\int_{0}^{\infty}M_{\beta}(\theta)S_{\alpha}(t^{\beta}\theta)d\theta,
\end{align*}
and
\begin{align*}
E_{\beta,\beta}(t)=\int_{0}^{\infty}\beta\theta M_{\beta}(\theta)S_{\alpha}(t^{\beta}\theta)d\theta,
\end{align*}
which contain the Mainardi's Wright-type function with $\beta\in (0,1)$ given by
\begin{align*}
M_{\beta}(\theta)=\sum_{n=0}^{\infty}\frac{(-1)^{n}\theta^{n}}{n!\Gamma(1-\beta(1+n))},
\end{align*}
in which the Mainardi function $M_{\beta}(\theta)$ act as a bridge between the classical integral-order and fractional derivatives of differential equations, for more details see [19,20]. Here, the derivation of mild solution (2.4) can be found in Appendix A.

Firstly, let us state the property of the special Mainardi function $M_{\beta}(\theta)$. Further, the properties of generalized Mittag-Leffler operators $E_{\beta}(t)$ and $E_{\beta,\beta}(t)$ are proved.

\textbf{Lemma 2.3.} (see [15]) For any $\beta\in (0,1)$ and $-1<\varepsilon<\infty$, it is not difficult to verity that
\begin{align*}
M_{\beta}(\theta)\geq0, ~and~ \int_{0}^{\infty}\theta^{\varepsilon}M_{\beta}(\theta)d\theta=\frac{\Gamma(1+\varepsilon)}{\Gamma(1+\beta\varepsilon)}, \tag{2.5}
\end{align*}
for all $\theta\geq0$.

\textbf{Theorem 2.1.} For any $t>0$, $E_{\beta}(t)$ and $E_{\beta,\beta}(t)$ are linear and bounded operators. Moreover, for $0\leq\nu<\alpha<2$, there exist constants $C_{\alpha}=C(\alpha,\beta,\nu)>0$ and $C_{\beta}=C(\alpha,\beta,\nu)>0$ such that
\begin{align*}
\|E_{\beta}(t)v\|_{\dot{H}^{\nu}}\leq C_{\alpha} t^{-\frac{\beta\nu}{\alpha}}\|v\|,\|E_{\beta,\beta}(t)v\|_{\dot{H}^{\nu}}\leq C_{\beta}t^{-\frac{\beta\nu}{\alpha}}\|v\|.\tag{2.6}
\end{align*}

\textbf{Proof.} For $t>0$ and $0\leq\nu<\alpha<2$, by means of the Lemma 2.1 and Lemma 2.3, we have
\begin{align*}
\|E_{\beta}(t)v\|_{\dot{H}^{\nu}}&\leq\int_{0}^{\infty}M_{\beta}(\theta)\|A_{\nu}S_{\alpha}(t^{\beta}\theta)v\|d\theta\\
&\leq\int_{0}^{\infty}C_{\alpha,\nu}t^{-\frac{\beta\nu}{\alpha}}\theta^{-\frac{\nu}{\alpha}}M_{\beta}(\theta) \|v\|d\theta\\
&=\frac{C_{\alpha,\nu}\Gamma(1-\frac{\nu}{\alpha})}{\Gamma(1-\frac{\beta\nu}{\alpha})}t^{-\frac{\beta\nu}{\alpha}}\|v\|, v\in L^{2}(D),
\end{align*}
and
\begin{align*}
\|E_{\beta,\beta}(t)v\|_{\dot{H}^{\nu}}&\leq\int_{0}^{\infty}\beta\theta M_{\beta}(\theta)\|A_{\nu}S_{\alpha}(t^{\beta}\theta)v\|d\theta\\
&\leq \int_{0}^{\infty}C_{\alpha,\nu}\beta t^{-\frac{\beta\nu}{\alpha}}\theta^{1-\frac{\nu}{\alpha}}M_{\beta}(\theta) \|v\|d\theta\\
&=\frac{C_{\alpha,\nu}\beta\Gamma(2-\frac{\nu}{\alpha})}{\Gamma(1+\beta(1-\frac{\nu}{\alpha}))}t^{-\frac{\beta\nu}{\alpha}}\|v\|, v\in L^{2}(D),
\end{align*}
which imply that the estimates (2.6) hold, so it is easy to know that $E_{\beta}(t)$ and $E_{\beta,\beta}(t)$ are linear and bounded operators.

\textbf{Theorem 2.2.} For any $t>0$, the operators $E_{\beta}(t)$ and $E_{\beta,\beta}(t)$ are strongly continuous. Moreover, for any $0\leq t_{1}< t_{2}\leq T$ and for $0<\nu<\alpha<2$, there exist constants $C_{\alpha\nu}=C(\alpha,\beta,\nu)>0$ and $C_{\beta\nu}=C(\alpha,\beta,\nu)>0$ such that
\begin{align*}
\|(E_{\beta}(t_{2})-E_{\beta}(t_{1}))v\|_{\dot{H}^{\nu}}\leq C_{\alpha\nu}(t_{2}-t_{1})^{\frac{\beta\nu}{\alpha}}\|v\|,\tag{2.7}
\end{align*}
and
\begin{align*}
\|(E_{\beta,\beta}(t_{2})-E_{\beta,\beta}(t_{1}))v\|_{\dot{H}^{\nu}}\leq C_{\beta\nu}(t_{2}-t_{1})^{\frac{\beta\nu}{\alpha}}\|v\|.\tag{2.8}
\end{align*}

\textbf{Proof.} For any $0<T_{0}\leq t_{1}<t_{2}\leq T$, it is easy to deduce that
\begin{align*}
\int_{t_{1}}^{t_{2}}\frac{d S_{\alpha}(t^{\beta}\theta)}{dt} dt&=S_{\alpha}(t_{2}^{\beta}\theta)-S_{\alpha}(t_{1}^{\beta}\theta)\\
&=-\int_{t_{1}}^{t_{2}} \beta t^{\beta-1}\theta A_{\alpha} S_{\alpha}(t^{\beta}\theta)dt.\tag{2.9}
\end{align*}

For $0<\nu<\alpha<2$, making use of the above expression, the Lemma 2.1 and Lemma 2.3, we can arrive at
\begin{align*}
\|(E_{\beta}(t_{2})-E_{\beta}(t_{1}))v\|_{\dot{H}^{\nu}}&=\|A_{\nu}(E_{\beta}(t_{2})-E_{\beta}(t_{1}))v\|\\
&=\|\int_{0}^{\infty}M_{\beta}(\theta)A_{\nu}(S_{\alpha}(t_{2}^{\beta}\theta)-S_{\alpha}(t_{1}^{\beta}\theta)vd\theta\|\\
&\leq\int_{0}^{\infty}\beta\theta M_{\beta}(\theta)\int_{t_{1}}^{t_{2}} t^{\beta-1}\| A_{\alpha+\nu} S_{\alpha}(t^{\beta}\theta)v\|_{L^{2}}dtd\theta\\
&\leq\int_{0}^{\infty}C_{\alpha,\nu}\beta\theta^{-\frac{\nu}{\alpha}}M_{\beta}(\theta)(\int_{t_{1}}^{t_{2}}t^{-\frac{\beta\nu}{\alpha}-1}dt) \|v\|d\theta\\
&=\frac{\alpha C_{\alpha,\nu}\Gamma(1-\frac{\nu}{\alpha})}{\nu\Gamma(1-\frac{\beta\nu}{\alpha})}(t_{1}^{-\frac{\beta\nu}{\alpha}}-t_{2}^{-\frac{\beta\nu}{\alpha}})\|v\|\\
&\leq \frac{\alpha C_{\alpha,\nu}\Gamma(1-\frac{\nu}{\alpha})}{\nu T_{0}^{\frac{2\beta\nu}{\alpha}}\Gamma(1-\frac{\beta\nu}{\alpha})}(t_{2}-t_{1})^{\frac{\beta\nu}{\alpha}}\|v\|, v\in L^{2}(D),
\end{align*}
and
\begin{align*}
\|(E_{\beta,\beta}(t_{2})-E_{\beta,\beta}(t_{1}))v\|_{\dot{H}^{\nu}}&=\|A_{\nu}(E_{\beta,\beta}(t_{2})-E_{\beta,\beta}(t_{1}))v\|\\
&=\|\int_{0}^{\infty}\beta\theta M_{\beta}(\theta)A_{\nu}(S_{\alpha}(t_{2}^{\beta}\theta)-S_{\alpha}(t_{1}^{\beta}\theta)vd\theta\|\\
&\leq\int_{0}^{\infty}\beta^{2}\theta^{2} M_{\beta}(\theta)\int_{t_{1}}^{t_{2}} t^{\beta-1}\| A_{\alpha+\nu} S_{\alpha}(t^{\beta}\theta)v\|dtd\theta\\
&\leq\int_{0}^{\infty}C_{\alpha,\nu}\beta^{2}\theta^{1-\frac{\nu}{\alpha}}M_{\beta}(\theta)(\int_{t_{1}}^{t_{2}}t^{-\frac{\beta\nu}{\alpha}-1}dt) \|v\|d\theta\\
&=\frac{\alpha\beta C_{\alpha,\nu}\Gamma(2-\frac{\nu}{\alpha})}{\nu\Gamma(1+\beta(1-\frac{\nu}{\alpha}))}(t_{1}^{-\frac{\beta\nu}{\alpha}}-t_{2}^{-\frac{\beta\nu}{\alpha}})\|v\|\\
&\leq \frac{\alpha\beta C_{\alpha,\nu}\Gamma(2-\frac{\nu}{\alpha})}{\nu T_{0}^{\frac{2\beta\nu}{\alpha}}\Gamma(1+\beta(1-\frac{\nu}{\alpha}))}(t_{2}-t_{1})^{\frac{\beta\nu}{\alpha}}\|v\|, v\in L^{2}(D).
\end{align*}

It is obviously to see that the term $\|(E_{\beta}(t_{2})-E_{\beta}(t_{1}))v\|_{\dot{H}^{\nu}}\rightarrow0$ and $\|(E_{\beta,\beta}(t_{2})-E_{\beta,\beta}(t_{1}))v\|_{\dot{H}^{\nu}}\rightarrow0$ as $t_{1}\rightarrow t_{2}$, which mean that the operators $E_{\beta}(t)$ and $E_{\beta,\beta}(t)$ are strongly continuous.

\textsl{Remark.} Assume $\nu=0$ in Theorem 2.2, then there exist constants $C_{\alpha}=C(\alpha,\beta)>0$ and $C_{\beta}=C(\alpha,\beta)>0$ such that
\begin{align*}
\|(E_{\beta}(t_{2})-E_{\beta}(t_{1}))v\|\leq C_{\alpha}(t_{2}-t_{1})\|v\|,\tag{2.10}
\end{align*}
and
\begin{align*}
\|(E_{\beta,\beta}(t_{2})-E_{\beta,\beta}(t_{1}))v\|\leq C_{\beta}(t_{2}-t_{1})\|v\|.\tag{2.11}
\end{align*}

\textbf{Proof.} For any $0<T_{0}\leq t_{1}<t_{2}\leq T$, the same as the proof of Theorem 2.2, we get
\begin{align*}
\|(E_{\beta}(t_{2})-E_{\beta}(t_{1}))v\|&=\|\int_{0}^{\infty}M_{\beta}(\theta)(S_{\alpha}(t_{2}^{\beta}\theta)-S_{\alpha}(t_{1}^{\beta}\theta)vd\theta\|_{L^{2}}\\
&\leq\int_{0}^{\infty}\beta\theta M_{\beta}(\theta)\int_{t_{1}}^{t_{2}} t^{\beta-1}\| A_{\alpha} S_{\alpha}(t^{\beta}\theta)v\|dtd\theta\\
&\leq \int_{0}^{\infty}C_{\alpha,\alpha}\beta M_{\beta}(\theta)(\int_{t_{1}}^{t_{2}}t^{-1}dt) \|v\|d\theta\\
&=C_{\alpha,\alpha}\beta(\ln t_{2}-\ln t_{1})\|v\|\\
&\leq \frac{C_{\alpha,\alpha}\beta}{T_{0}}(t_{2}-t_{1})\|v\|, v\in L^{2}(D),
\end{align*}
and
\begin{align*}
\|(E_{\beta,\beta}(t_{2})-E_{\beta,\beta}(t_{1}))v\|&=\|\int_{0}^{\infty}\beta\theta M_{\beta}(\theta)(S_{\alpha}(t_{2}^{\beta}\theta)-S_{\alpha}(t_{1}^{\beta}\theta)vd\theta\|\\
&\leq\int_{0}^{\infty}\beta^{2}\theta^{2} M_{\beta}(\theta)\int_{t_{1}}^{t_{2}} t^{\beta-1}\| A_{\alpha} S_{\alpha}(t^{\beta}\theta)v\|dtd\theta\\
&\leq\int_{0}^{\infty} C_{\alpha,\alpha}\beta^{2}\theta M_{\beta}(\theta)(\int_{t_{1}}^{t_{2}}t^{-1}dt) \|v\|d\theta\\
&=\frac{C_{\alpha,\alpha}\beta^{2}\Gamma(2)}{\Gamma(1+\beta)}(\ln t_{2}-\ln t_{1})\|v\|\\
&\leq \frac{C_{\alpha,\alpha}\beta^{2}\Gamma(2)}{T_{0}\Gamma(1+\beta)}(t_{2}-t_{1})\|v\|, v\in L^{2}(D).
\end{align*}

This completes the proof.

\section{Existence and uniqueness of mild solution}

Our main purpose of this section is to prove the existence and uniqueness of mild solution to the problem (2.1). To do this, the following assumptions are imposed.

\textbf{Assumption 3.1.} The measurable function $g: \Omega\times H \rightarrow L_{0}^{2}$ satisfies the following global Lipschitz and growth conditions:
\begin{align*}
\|g(v)\|_{L_{0}^{2}}\leq C\|v\|,\|g(u)-g(v)\|_{L_{0}^{2}}\leq C\|u-v\|,\tag{3.1}
\end{align*}
for any $u,v\in H$.

\textbf{Assumption 3.2.} Let $C>0$ be a real number, then the bounded bilinear operator $B:L^{2}(D)\rightarrow \dot{H}^{-1}(D)$ satisfies the following properties:
\begin{align*}
\|B(u)\|_{\dot{H}^{-1}}\leq C\|u\|^{2},\tag{3.2}
\end{align*}
and
\begin{align*}
\|B(u)-B(v)\|_{\dot{H}^{-1}}\leq C(\|u\|+\|v\|)\|u-v\|,\tag{3.3}
\end{align*}
for any $u,v\in L^{2}(D)$.

\textbf{Assumption 3.3.} Assume that the initial value $u_{0}: \Omega \rightarrow \dot{H}^{\nu}$ is a $\mathcal{F}_{0}-$measurable random variable, it holds that
\begin{align*}
\|u_{0}\|_{L^{p}(\Omega;\dot{H}^{\nu})}<\infty,  \tag{3.4}
\end{align*}
for any $0\leq\nu<\alpha<2$.

\textbf{Theorem 3.1.} Let Assumptions 3.1 to 3.3 be satisfied for some $p\geq2$, then there exists a unique mild solution $(u(t))_{t\in[0,T]}$ in the space $L^{p}(\Omega;\dot{H}^{\nu})$ with $0\leq\nu<\alpha<2$.

\textbf{Proof.} We fix an $\omega\in\Omega$ and use the standard Picard's iteration argument to prove the
existence of mild solution. To begin with, the sequence of stochastic process $\{u_{n}(t)\}_{n\geq0}$ is constructed as
\begin{align*}
\begin{cases}
u_{n+1}(t)=E_{\beta}(t)u_{0}+N_{1}(u_{n}(t))+N_{2}(u_{n}(t)),\\
u_{0}(t)=u_{0},
\end{cases} \tag{3.5}
\end{align*}
where
\begin{align*}
\begin{cases}
N_{1}(u_{n}(t))=\int_{0}^{t}(t-s)^{\beta-1}E_{\beta,\beta}(t-s)B(u_{n}(s))ds,\\
N_{2}(u_{n}(t))=\int_{0}^{t}(t-s)^{\beta-1}E_{\beta,\beta}(t-s)g(u_{n}(s))dW(s).
\end{cases} \tag{3.6}
\end{align*}

The proof will be split into three steps.

\textit{Step 1}: For each $n\geq0$, we show that
\begin{align*}
\sup\limits_{t\in[0,T]}\mathbb{E}[\|u_{n}(t)\|_{\dot{H}^{\nu}}^{p}]<\infty.
\end{align*}

Note that
\begin{align*}
\mathbb{E}[\|u_{n+1}(t)\|_{\dot{H}^{\nu}}^{p}]&\leq 3^{p-1} \mathbb{E}[\|E_{\beta}(t)u_{0}\|_{\dot{H}^{\nu}}^{p}]+3^{p-1}\mathbb{E}[\|N_{1}(u_{n}(t))\|_{\dot{H}^{\nu}}^{p}]\\
&+3^{p-1}\mathbb{E}[\|N_{2}(u_{n}(t))\|_{\dot{H}^{\nu}}^{p}].\tag{3.7}
\end{align*}

The application of the Lemma 2.1 gives
\begin{align*}
\mathbb{E}[\|E_{\beta}(t)u_{0}\|_{\dot{H}^{\nu}}]&\leq\mathbb{E}[\int_{0}^{\infty}M_{\beta}(\theta)(\|A_{\nu}S_{\alpha}(t^{\beta}\theta)u_{0}\|^{2})^{\frac{1}{2}}d\theta]\\
&=\mathbb{E}[\int_{0}^{\infty}M_{\beta}(\theta)(\sum\limits_{n=1}^{\infty}\langle A_{\nu}e^{-t^{\beta}\theta A_{\alpha}}u_{0},e_{n}\rangle^{2})^{\frac{1}{2}}d\theta]\\
&=\mathbb{E}[\int_{0}^{\infty}M_{\beta}(\theta)(\sum\limits_{n=1}^{\infty}\langle A_{\nu}u_{0},e^{-t^{\beta}\theta \lambda_{n}^{\frac{\alpha}{2}}}e_{n}\rangle^{2})^{\frac{1}{2}} d\theta]\\
&\leq \mathbb{E}[\int_{0}^{\infty}M_{\beta}(\theta) \|u_{0}\|_{\dot{H}^{\nu}} d\theta]=\mathbb{E}[\|u_{0}\|_{\dot{H}^{\nu}}].\tag{3.8}
\end{align*}

Applying the following H\"{o}lder inequality to the second term of the right-hand side of (3.7)
\begin{align*}
\int_{a}^{b}|f(t)g(t)|dt\leq(\int_{a}^{b}|f(t)|^{p}dt)^{\frac{1}{p}}(\int_{a}^{b}|g(t)|^{q}dt)^{\frac{1}{q}},\frac{1}{p}+\frac{1}{q}=1,
\end{align*}
where $p,q\in(1,\infty)$, we infer
\begin{align*}
\mathbb{E}[\|N_{1}(u_{n}(t))\|_{\dot{H}^{\nu}}^{p}]&\leq \mathbb{E}[(\int_{0}^{t}\|(t-s)^{\beta-1}A_{1}E_{\beta,\beta}(t-s)A_{\nu-1}B(u_{n}(s))\|ds)^{p}]\\
&\leq C_{\beta}^{p} (\int_{0}^{t}(t-s)^{\frac{p(\beta-1-\frac{\beta}{\alpha})}{p-1}}ds)^{p-1}\int_{0}^{t}\mathbb{E}[\|A_{\nu-1}B(u_{n}(s))\|^{p}]ds\\
&\leq K_{1}\int_{0}^{t}\mathbb{E}[\|u_{n}(s)\|_{\dot{H}^{\nu}}^{p}]ds, \tag{3.9}
\end{align*}
where $K_{1}=C_{\beta}^{p}C^{p}[\frac{p-1 }{p(\beta-\frac{\beta}{\alpha})-1}]^{p-1}T^{p(\beta-\frac{\beta}{\alpha})-1}(\max\limits_{t\in[0,T]}\mathbb{E}[\|u_{n}(t)\|_{\dot{H}^{\nu}}^{p}])$.

Making use of the H\"{o}lder inequality and Lemma 2.2 to the third term of the right-hand side of (3.7), we get
\begin{align*}
\mathbb{E}[\|N_{2}(u_{n}(t))\|_{\dot{H}^{\nu}}^{p}]&\leq C(p)\mathbb{E}[(\int_{0}^{t}\|(t-s)^{\beta-1}E_{\beta,\beta}(t-s)\|^{2}\|A_{\nu}g(u_{n}(s))\|_{L_{0}^{2}}^{2}ds)^{\frac{p}{2}}]\\
&\leq  C(p)C_{\beta}^{p} (\int_{0}^{t}(t-s)^{\frac{2p(\beta-1)}{p-2}}ds)^{\frac{p-2}{2}}\int_{0}^{t}\mathbb{E}\|A_{\nu}g(u_{n}(s))\|_{L_{0}^{2}}^{p}ds\\
&\leq K_{2}\int_{0}^{t}\mathbb{E}[\|u_{n}(s)\|_{\dot{H}^{\nu}}^{p}]ds, \tag{3.10}
\end{align*}
where $K_{2}=C(p)C_{\beta}^{p}C^{p}[\frac{p-2 }{p(2\beta-1)-2}]^{\frac{p-2}{2}}T^{\frac{p(2\beta-1)-2}{2}}$.

Using the above estimates (3.7)-(3.10), we have
\begin{align*}
\mathbb{E}[\|u_{n+1}(t)\|_{\dot{H}^{\nu}}^{p}]&\leq 3^{p-1}\mathbb{E}[\|u_{0}\|_{\dot{H}^{\nu}}^{p}]+3^{p-1}(K_{1}+K_{2})\int_{0}^{t}\mathbb{E}[\|u_{n}(s)\|_{\dot{H}^{\nu}}^{p}]ds.
\end{align*}

By means of the extension of Gronwall's lemma, it holds that
\begin{align*}
\sup\limits_{t\in[0,T]}\mathbb{E}[\|u_{n}(t)\|_{\dot{H}^{\nu}}^{p}]<\infty,
\end{align*}
for each $n\geq0$.

\textit{Step 2}: Show that the sequence $\{u_{n}(t)\}_{n\geq0}$ is a Cauchy sequence in the space $L^{p}(\Omega;\dot{H}^{\nu})$.

For any $n\geq m\geq1$, applying the similar arguments employed to obtain (3.9) and (3.10), we get
\begin{align*}
\mathbb{E}[\|u_{n}(t)-u_{m}(t)\|_{\dot{H}^{\nu}}^{p}]&\leq 2^{p-1}\mathbb{E}[\|N_{1}(u_{n-1}(t))-N_{1}(u_{m-1}(t))\|_{\dot{H}^{\nu}}^{p}]\\
&\hspace{2mm}+2^{p-1}\mathbb{E}[\|N_{2}(u_{n-1}(t))-N_{2}(u_{m-1}(t))\|_{\dot{H}^{\nu}}^{p}]\\
&\leq K\int_{0}^{t}\mathbb{E}[\|u_{n-1}(s)-u_{m-1}(s)\|_{\dot{H}^{\nu}}^{p}]ds,\tag{3.11}
\end{align*}
in which
\begin{align*}
K&=2^{p-1}\{C_{\beta}^{p}C^{p}[\frac{p-1 }{p(\beta-\frac{\beta}{\alpha})-1}]^{p-1}T^{p(\beta-\frac{\beta}{\alpha})-1}(\max\limits_{t\in[0,T]}\mathbb{E}[\|u_{n-1}(t)\|_{\dot{H}^{\nu}}^{p}]\\
&\hspace{2mm}+\max\limits_{t\in[0,T]}\mathbb{E}[\|u_{m-1}(t)\|_{\dot{H}^{\nu}}^{p}])+C(p)C_{\beta}^{p}C^{p}[\frac{p-2 }{p(2\beta-1)-2}]^{\frac{p-2}{2}}T^{\frac{p(2\beta-1)-2}{2}}\}.
\end{align*}

A direct application of Gronwall's lemma yields
\begin{align*}
\sup\limits_{t\in[0,T]}\mathbb{E}[\|u_{n}(t)-u_{m}(t)\|_{\dot{H}^{\nu}}^{p}]=0.
\end{align*}

As a result, the sequence $\{u_{n}(t)\}_{n\geq0}$ is a Cauchy sequence in the space $L^{p}(\Omega;\dot{H}^{\nu})$. Further, there exists a $u(t)\in L^{p}(\Omega;\dot{H}^{\nu})$ such that
\begin{align*}
\sup\limits_{t\in[0,T]}\mathbb{E}[\|u_{n}(t)-u(t)\|_{\dot{H}^{\nu}}^{p}]=0,
\end{align*}
for all $T>0$.

Taking limits to the stochastic sequence $\{u_{n}(t)\}_{n\geq0}$ in (3.5) as $n\rightarrow \infty$, we finish the proof of the existence of mild solution to (2.1).

\textit{Step 3}: We show the uniqueness of mild solution.

Assume $u$ and $v$ are two mild solutions of the problem (2.1), using the similar calculations as in Step 2, we can obtain
\begin{align*}
\sup\limits_{t\in[0,T]}\mathbb{E}[\|u(t)-v(t)\|_{\dot{H}^{\nu}}^{p}]=0,\tag{3.12}
\end{align*}
for all $T>0$, which implies that $u=v$, it follows that the uniqueness of mild solution.

Obviously, when $\nu=0$, the above three steps still work. Thus the proof of Theorem 3.1 is completed.

\section{Regularity of mild solution}

In this section, we will prove the spatial and temporal regularity properties of mild solution to time-space fractional SBE based on the analytic semigroup.

\textbf{Theorem 4.1.} Let Assumptions 3.1 to 3.3 hold with $1\leq\nu<\alpha<2$ and $p\geq2$, let $u(t)$ be a unique mild solution of the problem (2.1) with $\mathds{P}(u(t)\in\dot{H}^{\nu})=1$ for any $t\in [0,T]$, then there exists a constant $C$ such that
\begin{align*}
\sup\limits_{t\in[0,T]}\|u(t)\|_{L^{p}(\Omega;\dot{H}^{\nu})}\leq C(\|u_{0}\|_{L^{p}(\Omega;H)}+\sup\limits_{t\in[0,T]}\|u(t)\|_{L^{p}(\Omega;\dot{H}^{1})})\tag{4.1}
\end{align*}
\textbf{Proof.} For any $0\leq t\leq T$ and $1\leq\nu <\alpha< 2$, we have
\begin{align*}
\|u(t)\|_{L^{p}(\Omega;\dot{H}^{\nu})}&=(\mathbb{E}[\|u(t)\|_{\dot{H}^{\nu}}^{p}])^{\frac{1}{p}}=\|A_{\nu}u(t)\|_{L^{p}(\Omega;H)}\\
&\leq \|A_{\nu}E_{\beta}(t)u_{0}\|_{L^{p}(\Omega;H)}\\
&\hspace{2mm}+\|A_{\nu}\int_{0}^{t}(t-s)^{\beta-1}E_{\beta,\beta}(t-s)B(u(s))ds\|_{L^{p}(\Omega;H)}\\
&\hspace{2mm}+\|A_{\nu}\int_{0}^{t}(t-s)^{\beta-1}E_{\beta,\beta}(t-s)g(u(s))dW(s)\|_{L^{p}(\Omega;H)}\\
&=I+II+III. \tag{4.2}
\end{align*}

Using Theorem 2.1, the first term can be estimated by
\begin{align*}
I=\|A_{\nu}E_{\beta}(t)u_{0}\|_{L^{p}(\Omega;H)}\leq C_{\alpha} t^{-\frac{\beta\nu}{\alpha}}\|u_{0}\|_{L^{p}(\Omega;H)}<\infty.\tag{4.3}
\end{align*}

It is easy to know that
\begin{align*}
\int_{0}^{T} C_{\alpha} t^{-\frac{\beta\nu}{\alpha}}\|u_{0}\|_{L^{p}(\Omega;H)}dt=\frac{\alpha C_{\alpha}}{\alpha-\beta\nu}T^{\frac{\alpha-\beta\nu}{\alpha}}\|u_{0}\|_{L^{p}(\Omega;H)}.\tag{4.4}
\end{align*}

The application of Theorem 2.1 and Assumptions 3.2, we get
\begin{align*}
(II)^{p}&\leq\mathbb{E}[(\int_{0}^{t}\|(t-s)^{\beta-1}A_{\nu}E_{\beta,\beta}(t-s)B(u(s))\|ds)^{p}]\\
&\leq C_{\beta}^{p} (\int_{0}^{t}(t-s)^{\frac{p[\beta-1-\frac{\beta(\nu+1)}{\alpha}]}{p-1}}ds)^{p-1}\int_{0}^{t}\mathbb{E}[\|A_{-1}B(u(s))\|_{\dot{H}^{1}}^{p}]ds\\
&\leq C_{1}\sup\limits_{t\in[0,T]}\mathbb{E}[\|u(s)\|_{\dot{H}^{1}}^{p}], \tag{4.5}
\end{align*}
where $C_{1}=C_{\beta}^{p}C^{p}\{\frac{p-1 }{p[\beta-\frac{\beta(\nu+1)}{\alpha}]-1}\}^{p-1}T^{p[\beta-\frac{\beta(\nu+1)}{\alpha}]-1}(\max\limits_{t\in[0,T]}\mathbb{E}[\|u(t)\|_{\dot{H}^{1}}])$.

By means of Theorem 2.1, Assumptions 3.1 and Lemma 2.2, we can deduce
\begin{align*}
(III)^{p}&\leq C(p)\mathbb{E}[(\int_{0}^{t}\|(t-s)^{\beta-1}A_{\nu-1}E_{\beta,\beta}(t-s)\|^{2}\|A_{1}g(u(s))\|_{L_{2}^{0}}^{2}ds)^{\frac{p}{2}}]\\
&\leq  C(p)C_{\beta}^{p} (\int_{0}^{t}(t-s)^{\frac{2p[\beta-1-\frac{\beta(\nu-1)}{\alpha}]}{p-2}}ds)^{\frac{p-2}{2}}\int_{0}^{t}\mathbb{E}\|A_{1}g(u(s))\|_{L_{2}^{0}}^{p}ds\\
&\leq C_{2}\sup\limits_{t\in[0,T]}\mathbb{E}[\|u(s)\|_{\dot{H}^{1}}^{p}], \tag{4.6}
\end{align*}
where $C_{2}=C(p)C_{\beta}^{p}C^{p}[\frac{p-2 }{p(2\beta-1-\frac{\beta(\nu-1)}{\alpha})-2}]^{\frac{p-2}{2}}T^{\frac{p[2\beta-1-\frac{\beta(\nu-1)}{\alpha}]-2}{2}}$.

Thus, we conclude the proof of Theorem 4.1 by combining with the estimates (4.2)-(4.6).

Next, we will devote to the temporal regularity of the mild solution.

\textbf{Theorem 4.2.} Let Assumptions 3.1 to 3.3 be fulfilled with $0<\nu<\alpha<2$ and $p\geq2$, for any $0\leq t_{1}<t_{2}\leq T$, the unique mild solution $u(t)$ to the problem (2.1) is H\"{o}lder continuous with respect to the norm $\|\cdot\|_{L^{p}(\Omega;\dot{H}^{\nu})}$ and satisfies
\begin{align*}
\|u(t_{2})-u(t_{1})\|_{L^{p}(\Omega;\dot{H}^{\nu})}\leq C(t_{2}-t_{1})^{\gamma}.\tag{4.7}
\end{align*}
\textbf{Proof.} For any $0\leq t_{1}<t_{2}\leq T$, from the mild solution (2.4), we have
\begin{align*}
u(t_{2})-u(t_{1})&=E_{\beta}(t_{2})u_{0}-E_{\beta}(t_{1})u_{0}\\
&\hspace{2mm}+\int_{0}^{t_{2}}(t_{2}-s)^{\beta-1}E_{\beta,\beta}(t_{2}-s)B(u(s))ds\\
&\hspace{2mm}-\int_{0}^{t_{1}}(t_{1}-s)^{\beta-1}E_{\beta,\beta}(t_{1}-s)B(u(s))ds\\
&\hspace{2mm}+\int_{0}^{t_{2}}(t_{2}-s)^{\beta-1}E_{\beta,\beta}(t_{2}-s)g(u(s))dW(s)\\
&\hspace{2mm}-\int_{0}^{t_{1}}(t_{1}-s)^{\beta-1}E_{\beta,\beta}(t_{1}-s)g(u(s))dW(s)\\
&=: I_{1}+I_{2}+I_{3},\tag{4.8}
\end{align*}
where
\begin{align*}
I_{1}=E_{\beta}(t_{2})u_{0}-E_{\beta}(t_{1})u_{0},
\end{align*}
\begin{align*}
I_{2}&=\int_{0}^{t_{2}}(t_{2}-s)^{\beta-1}E_{\beta,\beta}(t_{2}-s)B(u(s))ds\\
&\hspace{2mm}-\int_{0}^{t_{1}}(t_{1}-s)^{\beta-1}E_{\beta,\beta}(t_{1}-s)B(u(s))ds\\
&=\int_{0}^{t_{1}}(t_{1}-s)^{\beta-1}[E_{\beta,\beta}(t_{2}-s)-E_{\beta,\beta}(t_{1}-s)]B(u(s))ds\\
&\hspace{2mm}+\int_{0}^{t_{1}}[(t_{2}-s)^{\beta-1}-(t_{1}-s)^{\beta-1}]E_{\beta,\beta}(t_{2}-s)B(u(s))ds\\
&\hspace{2mm}+\int_{t_{1}}^{t_{2}}(t_{2}-s)^{\beta-1}E_{\beta,\beta}(t_{2}-s)B(u(s))ds\\
&=: I_{21}+I_{22}+I_{23}.\tag{4.9}
\end{align*}
and
\begin{align*}
I_{3}&=\int_{0}^{t_{2}}(t_{2}-s)^{\beta-1}E_{\beta,\beta}(t_{2}-s)g(u(s))dW(s)\\
&\hspace{2mm}-\int_{0}^{t_{1}}(t_{1}-s)^{\beta-1}E_{\beta,\beta}(t_{1}-s)g(u(s))dW(s)\\
&=\int_{0}^{t_{1}}(t_{1}-s)^{\beta-1}[E_{\beta,\beta}(t_{2}-s)-E_{\beta,\beta}(t_{1}-s)]g(u(s))dW(s)\\
&\hspace{2mm}+\int_{0}^{t_{1}}[(t_{2}-s)^{\beta-1}-(t_{1}-s)^{\beta-1}]E_{\beta,\beta}(t_{2}-s)g(u(s))dW(s)\\
&\hspace{2mm}+\int_{t_{1}}^{t_{2}}(t_{2}-s)^{\beta-1}E_{\beta,\beta}(t_{2}-s)g(u(s))dW(s)\\
&=: I_{31}+I_{32}+I_{33}.\tag{4.10}
\end{align*}

For any $0<\nu<\alpha<2$ and $p\geq2$, by virtue of Theorem 2.2, it follows that
\begin{align*}
\mathbb{E}[\|I_{1}\|_{\dot{H}^{\nu}}^{p}]&=\mathbb{E}[\|A_{\nu}[E_{\beta}(t_{2})-E_{\beta}(t_{1})]u_{0}\|^{p}] \\
&\leq C_{\alpha\nu}^{p}(t_{2}-t_{1})^{\frac{p\beta\nu}{\alpha}}\mathbb{E}[\|u_{0}\|^{p}].\tag{4.11}
\end{align*}

For the first term $I_{21}$ in (4.9), applying the Assumptions 3.2 and Theorem 2.2 and H\"{o}lder's inequality, we have
\begin{align*}
\mathbb{E}[\|I_{21}\|_{\dot{H}^{\nu}}^{p}]&=\mathbb{E}[\|\int_{0}^{t_{1}}(t_{1}-s)^{\beta-1}A_{\nu}[E_{\beta,\beta}(t_{2}-s)-E_{\beta,\beta}(t_{1}-s)]B(u(s))ds\|^{p}]\\
&\leq C_{\beta\nu}^{p}(t_{2}-t_{1})^{\frac{p\beta(\nu+1)}{\alpha}}(\int_{0}^{t_{1}}(t_{1}-s)^{\frac{p(\beta-1)}{p-1}}ds)^{p-1}\int_{0}^{t}\mathbb{E}[\|A_{-1}B(u(s))\|_{\dot{H}^{1}}^{p}]ds\\
&\leq C^{p}C_{\beta\nu}^{p}T^{p\beta}(\frac{p-1}{p\beta-1})^{p-1}(\sup\limits_{t\in[0,T]}\mathbb{E}[\|u(s)\|_{\dot{H}^{1}}^{2p}])(t_{2}-t_{1})^{\frac{p\beta(\nu+1)}{\alpha}}.\tag{4.12}
\end{align*}

Using the Assumptions 3.2, Theorem 2.1 and H\"{o}lder's inequality, we get
\begin{align*}
\mathbb{E}[\|I_{22}\|_{\dot{H}^{\nu}}^{p}]&=\mathbb{E}[\|\int_{0}^{t_{1}}[(t_{2}-s)^{\beta-1}-(t_{1}-s)^{\beta-1}]A_{\nu}E_{\beta,\beta}(t_{2}-s)B(u(s))ds\|^{p}]\\
&\leq C_{\beta}^{p}(\int_{0}^{t_{1}}\{[(t_{2}-s)^{\beta-1}-(t_{1}-s)^{\beta-1}]\times(t_{2}-s)^{-\frac{\beta(\nu+1)}{\alpha}}\}^{\frac{p}{p-1}}ds)^{p-1}\\
&\hspace{2mm}\times\int_{0}^{t_{1}}\mathbb{E}[\|A_{-1}B(u(s))\|_{\dot{H}^{1}}^{p}]ds\\
&\leq C^{p}C_{\beta}^{p}T\{\frac{p-1 }{p[\beta-\frac{\beta(\nu+1)}{\alpha}]-1}\}^{p-1}(\sup\limits_{t\in[0,T]}\mathbb{E}[\|u(s)\|_{\dot{H}^{1}}^{2p}])(t_{2}-t_{1})^{\frac{p\beta(\alpha-\nu-1)-\alpha}{\alpha}},\tag{4.13}
\end{align*}
and
\begin{align*}
\mathbb{E}[\|I_{23}\|_{\dot{H}^{\nu}}^{p}]&=\mathbb{E}[\|\int_{t_{1}}^{t_{2}}(t_{2}-s)^{\beta-1}A_{\nu}E_{\beta,\beta}(t_{2}-s)B(u(s))ds\|^{p}]\\
&\leq C_{\beta}^{p}(\int_{t_{1}}^{t_{2}}[(t_{2}-s)^{\beta-1-\frac{\beta(\nu+1)}{\alpha}}]^{\frac{p}{p-1}}ds)^{p-1}\int_{t_{1}}^{t_{2}}\mathbb{E}[\|A_{-1}B(u(s))\|_{\dot{H}^{1}}^{p}]ds\\
&\leq C^{p}C_{\beta}^{p}\{\frac{p-1 }{p[\beta-\frac{\beta(\nu+1)}{\alpha}]-1}\}^{p-1}(\sup\limits_{t\in[0,T]}\mathbb{E}[\|u(s)\|_{\dot{H}^{1}}^{2p}])(t_{2}-t_{1})^{\frac{p\beta(\alpha-\nu-1)}{\alpha}}.\tag{4.14}
\end{align*}

Next, by following the similar arguments as in the proof of (4.12)-(4.14) and using the Lemma 2.2, there holds
\begin{align*}
\mathbb{E}[\|I_{31}\|_{\dot{H}^{\nu}}^{p}]&=\mathbb{E}[\|\int_{0}^{t_{1}}(t_{1}-s)^{\beta-1}A_{\nu}[E_{\beta,\beta}(t_{2}-s)-E_{\beta,\beta}(t_{1}-s)]g(u(s))dW(s)\|^{p}]\\
&\leq C(p)\mathbb{E}[(\int_{0}^{t_{1}}\|(t_{1}-s)^{\beta-1}A_{\nu}[E_{\beta,\beta}(t_{2}-s)-E_{\beta,\beta}(t_{1}-s)]\|^{2}\|g(u(s))\|_{L_{0}^{2}}^{2}ds)^{\frac{p}{2}}]\\
&\leq  C(p)C_{\beta\nu}^{p}(t_{2}-t_{1})^{\frac{p\beta\nu}{\alpha}}(\int_{0}^{t_{1}}(t_{1}-s)^{\frac{2p(\beta-1)}{p-2}}ds)^{\frac{p-2}{2}}\int_{0}^{t_{1}}\mathbb{E}\|g(u(s))\|_{L_{0}^{2}}^{p}ds\\
&\leq C(p)C^{p}C_{\beta\nu}^{p}T^{\frac{2p\beta-p-1}{2}}(\frac{p-1}{2p\beta-p-2})^{p-1}(\sup\limits_{t\in[0,T]}\mathbb{E}[\|u(t)\|^{p}])(t_{2}-t_{1})^{\frac{p\beta\nu}{\alpha}},\tag{4.15}
\end{align*}
and
\begin{align*}
\mathbb{E}[\|I_{32}\|_{\dot{H}^{\nu}}^{p}]&=\mathbb{E}[\|\int_{0}^{t_{1}}[(t_{2}-s)^{\beta-1}-(t_{1}-s)^{\beta-1}]A_{\nu}E_{\beta,\beta}(t_{2}-s)g(u(s))dW(s)\|^{p}]\\
&\leq C(p)\mathbb{E}[(\int_{0}^{t_{1}}\|[(t_{2}-s)^{\beta-1}-(t_{1}-s)^{\beta-1}]A_{\nu}E_{\beta,\beta}(t_{2}-s)\|^{2}\|g(u(s))\|_{L_{0}^{2}}^{2}ds)^{\frac{p}{2}}]\\
&\leq C(p)C_{\beta}^{p}(\int_{0}^{t_{1}}\{[(t_{2}-s)^{\beta-1}-(t_{1}-s)^{\beta-1}]\times(t_{2}-s)^{-\frac{\beta\nu}{\alpha}}\}^{\frac{2p}{p-2}}ds)^{\frac{p-2}{2}}\\
&\hspace{2mm}\times\int_{0}^{t_{1}}\mathbb{E}\|g(u(s))\|_{L_{0}^{2}}^{p}ds\\
&\leq C(p)C^{p}C_{\beta}^{p}T[\frac{\alpha(p-2)}{2p\beta(\alpha-\nu)-(p+2)\alpha}]^{\frac{p-2}{2}}(\sup\limits_{t\in[0,T]}\mathbb{E}[\|u(t)\|^{p}])(t_{2}-t_{1})^{\frac{2p\beta(\alpha-\nu)-(p+2)\alpha}{2\alpha}},\tag{4.16}
\end{align*}
and
\begin{align*}
\mathbb{E}[\|I_{33}\|_{\dot{H}^{\nu}}^{p}]&=\mathbb{E}[\|\int_{t_{1}}^{t_{2}}(t_{2}-s)^{\beta-1}A_{\nu}E_{\beta,\beta}(t_{2}-s)g(u(s))dW(s)\|^{p}]\\
&\leq C(p)\mathbb{E}[(\int_{t_{1}}^{t_{2}}\|(t_{2}-s)^{\beta-1}A_{\nu}E_{\beta,\beta}(t_{2}-s)\|^{2}\|g(u(s))\|_{L_{0}^{2}}^{2}ds)^{\frac{p}{2}}]\\
&\leq C(p)C_{\beta}^{p}(\int_{t_{1}}^{t_{2}}[(t_{2}-s)^{\beta-1-\frac{\beta\nu}{\alpha}}]^{\frac{2p}{p-2}}ds)^{\frac{p-2}{2}}\int_{t_{1}}^{t_{2}}\mathbb{E}\|g(u(s))\|_{L_{0}^{2}}^{p}ds\\
&\leq C(p)C^{p}C_{\beta}^{p}[\frac{\alpha(p-2)}{2p\beta(\alpha-\nu)-(p+2)\alpha}]^{\frac{p-2}{2}}(\sup\limits_{t\in[0,T]}\mathbb{E}[\|u(t)\|^{p}])(t_{2}-t_{1})^{\frac{2p\beta(\alpha-\nu)-p\alpha}{2\alpha}}.\tag{4.17}
\end{align*}

Taking expectation on the both side of (4.8), and in view of the estimates (4.11)-(4.17), we conclude that
\begin{align*}
\|u(t_{2})-u(t_{1})\|_{L^{p}(\Omega;\dot{H}^{\nu})}\leq C(t_{2}-t_{1})^{\gamma},\tag{4.18}
\end{align*}
in which we take $\gamma=\min\{\frac{\beta\nu}{\alpha},\frac{p\beta(\alpha-\nu-1)-\alpha}{p\alpha},\frac{2p\beta(\alpha-\nu)-(p+2)\alpha}{2p\alpha}\}$ when $0<t_{2}-t_{1}<1$. Otherwise, if $t_{2}-t_{1}\geq1$, then we set $\gamma=\max\{\frac{\beta(\nu+1)}{\alpha},\frac{\beta(\alpha-\nu-1)}{\alpha},\frac{2p\beta(\alpha-\nu)-p\alpha}{2p\alpha}\}$.  

This completes the proof of Theorem 4.2.

\section*{Acknowledgements}

The authors are very thankful to the anonymous referees for their valuable comments and constructive suggestions, which helped us to improve the manuscript.
 Guang-an Zou is supported by National Nature Science Foundation of China (Grant No. 11626085), Bo Wang is supported by the foundation for Young
University Key Teacher by the Educational Department of Henan Province (No.2014GGJS-021).

\newpage
\section*{Appendix A}

Considering the following abstract formulation of time-space fractional stochastic Burgers equation:
\begin{align*}
\begin{cases}
^{C}D_{t}^{\beta}u(t)=-A_{\alpha}u(t)+B(u(t))+g(u(t))\frac{d W(t)}{d t},t>0,\\
u(0)=u_{0}.
\end{cases} \tag{A1}
\end{align*}

We derive the mild solution to (A1) by means of Laplace transform, which denoted by ~$\widehat{}$~.
Let $\lambda>0$, and we define that
\begin{align*}
\widehat{u}(\lambda)=\int_{0}^{\infty}e^{-\lambda s}u(s)ds,~ \widehat{B}(\lambda)=\int_{0}^{\infty}e^{-\lambda s}B(u(s))ds,
\end{align*}
and
\begin{align*}
\widehat{G}(\lambda)=\int_{0}^{\infty}e^{-\lambda s}[g(u(s))\frac{d W(s)}{d s}]ds=\int_{0}^{\infty}e^{-\lambda s}g(u(s))d W(s).
\end{align*}

Upon Laplace transform, using the formula $\widehat{^{C}D_{t}^{\beta}u}=\lambda^{\beta}\widehat{u}-\lambda^{\beta-1}u_{0}$. Then applying the Laplace transform to (A1), we obtain
\begin{align*}
\widehat{u}(\lambda)&=\frac{1}{\lambda}u_{0}+\frac{1}{\lambda^{\beta}}(-A_{\alpha})\widehat{u}(\lambda)+\frac{1}{\lambda^{\beta}}[\widehat{B}(\lambda)+\widehat{G}(\lambda)]\\
&=\lambda^{\beta-1}(\lambda^{\beta}I+A_{\alpha})^{-1}u_{0}+(\lambda^{\beta}I+A_{\alpha})^{-1}[\widehat{B}(\lambda)+\widehat{G}(\lambda)]\\
&=\lambda^{\beta-1}\int_{0}^{\infty}e^{-\lambda^{\beta}s}S_{\alpha}(s)u_{0}ds+\int_{0}^{\infty}e^{-\lambda^{\beta}s}S_{\alpha}(s)[\widehat{B}(\lambda)+\widehat{G}(\lambda)]ds,\tag{A2}
\end{align*}
in which $I$ is the identity operator, and $S_{\alpha}(t)=e^{-tA_{\alpha}}$ is an analytic semigroup generated by the operator $-A_{\alpha}$.

We introduce the following one-sided stable probability density function:
\begin{align*}
\omega_{\beta}(\theta)=\frac{1}{\pi}\sum_{n=1}^{\infty}(-1)^{n-1}\theta^{-\beta n-1}\frac{\Gamma(\beta n+1)}{n!}\sin(n\pi \beta), \theta\in(0,+\infty),\tag{A3}
\end{align*}
whose Laplace transform is given by
\begin{align*}
\int_{0}^{\infty}e^{-\lambda \theta}\omega_{\beta}(\theta)d\theta=e^{-\lambda^{\beta}}, 0<\beta<1.\tag{A4}
\end{align*}

Making use of above expression (A4), then the terms on the right-hand side of (A2) can be written as
\begin{align*}
&\lambda^{\beta-1}\int_{0}^{\infty}e^{-\lambda^{\beta}s}S_{\alpha}(s)u_{0}ds\\
&=\int_{0}^{\infty}\lambda^{\beta-1}e^{-\lambda^{\beta}t^{\beta}}S_{\alpha}(t^{\beta})u_{0}d(t^{\beta})\\
&=\int_{0}^{\infty}\beta(\lambda t)^{\beta-1}e^{-(\lambda t)^{\beta}}S_{\alpha}(t^{\beta})u_{0}dt\\
&=\int_{0}^{\infty}-\frac{1}{\lambda}\frac{d}{dt}[e^{-(\lambda t)^{\beta}}]S_{\alpha}(t^{\beta})u_{0}dt\\
&=\int_{0}^{\infty}\int_{0}^{\infty}\theta\omega_{\beta}(\theta)e^{-\lambda t\theta}S_{\alpha}(t^{\beta})u_{0}d\theta dt\\
&=\int_{0}^{\infty}e^{-\lambda t}[\int_{0}^{\infty}\omega_{\beta}(\theta)S_{\alpha}(\frac{t^{\beta}}{\theta^{\beta}})u_{0}d\theta ]dt,\tag{A5}
\end{align*}
and
\begin{align*}
&\int_{0}^{\infty}e^{-\lambda^{\beta}s}S_{\alpha}(s)\widehat{B}(\lambda)ds\\
&=\int_{0}^{\infty}\beta t^{\beta-1}e^{-(\lambda t)^{\beta}}S_{\alpha}(t^{\beta})\widehat{B}(\lambda)dt\\
&=\int_{0}^{\infty}\int_{0}^{\infty}\beta t^{\beta-1}e^{-(\lambda t)^{\beta}}S_{\alpha}(t^{\beta})e^{-\lambda s}B(u(s))dsdt\\
&=\int_{0}^{\infty}\int_{0}^{\infty}\int_{0}^{\infty}\beta\omega_{\beta}(\theta)e^{-\lambda t\theta}S_{\alpha}(t^{\beta})e^{-\lambda s}t^{\beta-1}B(u(s))d\theta dsdt\\
&=\int_{0}^{\infty}\int_{0}^{\infty}\int_{0}^{\infty}\beta\omega_{\beta}(\theta)e^{-\lambda(t+s)}S_{\alpha}(\frac{t^{\beta}}{\theta^{\beta}})\frac{t^{\beta-1}}{\theta^{\beta}}B(u(s))d\theta dsdt\\
&=\int_{0}^{\infty}e^{-\lambda t}[\beta \int_{0}^{t}\int_{0}^{\infty}\omega_{\beta}(\theta)S_{\alpha}(\frac{(t-s)^{\beta}}{\theta^{\beta}})\frac{(t-s)^{\beta-1}}{\theta^{\beta}}B(u(s))d\theta d s]dt.\tag{A6}
\end{align*}
and
\begin{align*}
&\int_{0}^{\infty}e^{-\lambda^{\beta}s}S_{\alpha}(s)\widehat{G}(\lambda)ds\\
&=\int_{0}^{\infty}\beta t^{\beta-1}e^{-(\lambda t)^{\beta}}S_{\alpha}(t^{\beta})\widehat{G}(\lambda)dt\\
&=\int_{0}^{\infty}\int_{0}^{\infty}\beta t^{\beta-1}e^{-(\lambda t)^{\beta}}S_{\alpha}(t^{\beta})e^{-\lambda s}g(u(s))d W(s)dt\\
&=\int_{0}^{\infty}\int_{0}^{\infty}\int_{0}^{\infty}\beta\omega_{\beta}(\theta)e^{-\lambda t\theta}S_{\alpha}(t^{\beta})e^{-\lambda s}t^{\beta-1}g(u(s))d\theta d W(s)dt\\
&=\int_{0}^{\infty}\int_{0}^{\infty}\int_{0}^{\infty}\beta\omega_{\beta}(\theta)e^{-\lambda(t+s)}S_{\alpha}(\frac{t^{\beta}}{\theta^{\beta}})\frac{t^{\beta-1}}{\theta^{\beta}}g(u(s))d\theta d W(s)dt\\
&=\int_{0}^{\infty}e^{-\lambda t}[\beta \int_{0}^{t}\int_{0}^{\infty}\omega_{\beta}(\theta)S_{\alpha}(\frac{(t-s)^{\beta}}{\theta^{\beta}})\frac{(t-s)^{\beta-1}}{\theta^{\beta}}g(u(s))d\theta d W(s)]dt.\tag{A7}
\end{align*}

Together with (A2) and (A5)-(A7) helps us to get
\begin{align*}
\widehat{u}(\lambda)&=\int_{0}^{\infty}e^{-\lambda t}[\int_{0}^{\infty}\omega_{\beta}(\theta)S_{\alpha}(\frac{t^{\beta}}{\theta^{\beta}})u_{0}d\theta ]dt\\
&\hspace{2mm}+\int_{0}^{\infty}e^{-\lambda t}[\beta \int_{0}^{t}\int_{0}^{\infty}\omega_{\beta}(\theta)S_{\alpha}(\frac{(t-s)^{\beta}}{\theta^{\beta}})\frac{(t-s)^{\beta-1}}{\theta^{\beta}}B(u(s))d\theta d s]dt\\
&\hspace{2mm}+\int_{0}^{\infty}e^{-\lambda t}[\beta \int_{0}^{t}\int_{0}^{\infty}\omega_{\beta}(\theta)S_{\alpha}(\frac{(t-s)^{\beta}}{\theta^{\beta}})\frac{(t-s)^{\beta-1}}{\theta^{\beta}}g(u(s))d\theta d W(s)]dt.\tag{A8}
\end{align*}

Now, by means of inverse Laplace transform to (A8), we have achieved that
\begin{align*}
u(t)&=\int_{0}^{\infty}\omega_{\beta}(\theta)S_{\alpha}(\frac{t^{\beta}}{\theta^{\beta}})u_{0}d\theta\\
&\hspace{2mm}+\beta \int_{0}^{t}\int_{0}^{\infty}\omega_{\beta}(\theta)S_{\alpha}(\frac{(t-s)^{\beta}}{\theta^{\beta}})\frac{(t-s)^{\beta-1}}{\theta^{\beta}}B(u(s))d\theta d s\\
&\hspace{2mm}+\beta \int_{0}^{t}\int_{0}^{\infty}\omega_{\beta}(\theta)S_{\alpha}(\frac{(t-s)^{\beta}}{\theta^{\beta}})\frac{(t-s)^{\beta-1}}{\theta^{\beta}}g(u(s))d\theta d W(s)\\
&=\int_{0}^{\infty}\frac{1}{\beta}\theta^{-\frac{1}{\beta}-1}\omega_{\beta}(\theta^{-\frac{1}{\beta}})S_{\alpha}(t^{\beta}\theta)u_{0}d\theta\\
&\hspace{2mm}+\int_{0}^{t}\int_{0}^{\infty}\theta^{-\frac{1}{\beta}}\omega_{\beta}(\theta^{-\frac{1}{\beta}})S_{\alpha}((t-s)^{\beta}\theta)(t-s)^{\beta-1}B(u(s))d\theta d s\\
&\hspace{2mm}+\int_{0}^{t}\int_{0}^{\infty}\theta^{-\frac{1}{\beta}}\omega_{\beta}(\theta^{-\frac{1}{\beta}})S_{\alpha}((t-s)^{\beta}\theta)(t-s)^{\beta-1}g(u(s))d\theta d W(s). \tag{A9}
\end{align*}

Here, we also introduce the Mainardi's Wright-type function
\begin{align*}
M_{\beta}(\theta)&=\sum_{n=0}^{\infty}\frac{(-1)^{n}\theta^{n}}{n!\Gamma(1-\beta(1+n))}\\
&=\frac{1}{\pi}\sum_{n=1}^{\infty}\frac{(-1)^{n-1}\theta^{n-1}}{(n-1)!}\Gamma(n\beta)\sin(n\pi \beta),
\end{align*}
where $0<\beta<1$ and $\theta\in(0,+\infty)$. Further, the relationships between the probability density function $\omega_{\beta}(\theta)$ and Mainardi's Wright-type function $M_{\beta}(\theta)$ are shown that
\begin{align*}
M_{\beta}(\theta)=\frac{1}{\beta}\theta^{-\frac{1}{\beta}-1}\omega_{\beta}(\theta^{-\frac{1}{\beta}}).
\end{align*}

We denote the generalized Mittag-Leffler operators $E_{\beta}(t)$ and $E_{\beta,\beta}(t)$ as
\begin{align*}
E_{\beta}(t)=\int_{0}^{\infty}M_{\beta}(\theta)S_{\alpha}(t^{\beta}\theta)d\theta,
\end{align*}
and
\begin{align*}
E_{\beta,\beta}(t)=\int_{0}^{\infty}\beta\theta M_{\beta}(\theta)S_{\alpha}(t^{\beta}\theta)d\theta.
\end{align*}

Therefore, the equation (A9) can be written as
\begin{align*}
 u(t)&=E_{\beta}(t)u_{0}+\int_{0}^{t}(t-s)^{\beta-1}E_{\beta,\beta}(t-s)B(u(s))ds\\
 &\hspace{2mm}+\int_{0}^{t}(t-s)^{\beta-1}E_{\beta,\beta}(t-s)g(u(s))dW(s). \tag{A10}
\end{align*}

Up to now, we have deduced the mild solution (A10) to the time-space fractional stochastic Burgers equation (A1).

\newpage

\section*{References}

[1] C. Gugg, H. Kielh\"{o}fer, M. Niggemann, On the approximation of the stochastic Burgers equation, Commun. Math. Phys. 230(1) (2002) 181-199.

[2] Z. Dong, T.G. Xu, One-dimensional stochastic Burgers equation driven by L\'{e}vy processes, J. Funct. Anal. 243(2) (2007) 631-678.

[3] Z. Dong, L. Xu, X. Zhang, Exponential ergodicity of stochastic Burgers equations driven by $\alpha$-stable processes, J. Stat. Phys. 154(4) (2014) 929-949.

[4] E. Weinan, K. Khanin, A. Mazel, et al, Invariant measure for Burgers equation with stochastic forcing, Ann. Math. 151(3) (2000) 877-960.

[5] N. Sugimoto, Generalized Burgers equations and fractional calculus,
Nonlinear Wave Motion.(A. Jeffery, Ed), (1991) 162-179.

[6] Z. Brze\'{z}niak, L. Debbi, B. Goldys, Ergodic properties of fractional stochastic Burgers equation, arXiv preprint arXiv:1106.1918, 2011.

[7] M. Inc, The approximate and exact solutions of the space- and time-fractional Burgers equations with initial conditions by variational iteration method, J. Math. Anal. Appl. 345(1) (2008) 476-484.

[8] S. Momani, Non-perturbative analytical solutions of the space- and time-fractional Burgers equations, Chaos Soliton Fract. 28 (2006) 930-937.

[9] H.M. Srivastava, J.J. Trujillo, Theory and applications of fractional differential equations, Elsevier, 2006.

[10] Z. Brze\'{z}niak, L. Debbi, On stochastic Burgers equation driven by a fractional Laplacian and space-time white noise, Stochastic differential equations: theorem and applications Interdiscip Math. Sci., vol. 2, World Sci. Publ., Hackensack, NJ (2007) 135-167. 

[11] D. Yang, m-Dissipativity for Kolmogorov operator of a fractional Burgers equation with space-time white noise, Potential Anal. 44(2) (2016) 215-227.

[12] G. Lv, J. Duan, Martingale and weak solutions for a stochastic nonlocal Burgers equation on finite intervals, J. Math. Anal. Appl. 449(1) (2017) 176-194.

[13] R. Kruse, Strong and weak approximation of semilinear stochastic evolution equations, Springer, 2014.

[14] R.N. Wang, D.H. Chen, T.J. Xiao, Abstract fractional Cauchy problems with almost sectorial operators, J. Differential Equations 252(1) (2012) 202-235.

[15] P.M. De Carvalho-Neto, P. Gabriela, Mild solutions to the time fractional Navier-Stokes equations in $R^{N}$, J. Differential Equations 259 (2015) 2948-2980.

[16] Y. Zhou, L. Peng, On the time-fractional Navier-Stokes equations. Comput. Math. Appl. 73(6) (2017) 874-891.

[17] R. Sakthivel, S. Suganya, S.M. Anthoni, Approximate controllability of fractional stochastic evolution equations, Comput. Math. Appl. 63 (2012) 660-668.

[18] Y. Zhou, F. Jiao, Existence of mild solutions for fractional neutral evolution equations, Comput. Math. Appl. 59 (2010) 1063-1077.

[19] F. Mainardi, On the initial value problem for the fractional diffusion-wave equation. Waves and Stability in Continuous Media, World Scientific, Singapore, 1994.

[20] F. Mainardi, The fundamental solutions for the fractional diffusion-wave equation, Appl. Math. Lett. 9(6) (1996) 23-28.

[21] X.J. Yang, H.M. Srivastava, J.A. Machado, A new fractional derivative without singular kernel: Application to the modelling of the steady heat flow, Therm. Sci. 20(2) (2016) 753-756.

[22] X.J. Yang, J.A.T. Machado, J. Hristov, Nonlinear dynamics for local fractional Burgers¡¯equation arising in fractal flow, Nonlinear Dynam. 84(1) (2016) 3-7.

[23] X.J. Yang, F. Gao, H.M. Srivastava, Exact travelling wave solutions for the local fractional two-dimensional Burgers-type equations, Comput. Math. Appl. 73(2) (2017) 203-210.

[24] X.J. Yang, Advanced local fractional calculus and its applications, World Science, New York, 2012.

[25] L. Debbi, Well-posedness of the multidimensional fractional stochastic Navier-Stokes equations on the torus and on bounded domains, J. Math. Fluid Mech. 18(1) (2016) 25-69.

\end{document}